\documentclass[10pt,a4paper]{amsart}
\usepackage{amssymb}
\usepackage[english]{babel}
\usepackage{epsfig}
\usepackage{amsthm}
\usepackage{amsmath}
\usepackage{amssymb}
\usepackage{wasysym}
\usepackage{tikz}
\usepackage[all]{xy}
\usepackage{url}
\newtheorem{teo}{Theorem}[section]
\newtheorem{pro}[teo]{Proposition}
\newtheorem{lem}[teo]{Lemma}
\newtheorem{cor}[teo]{Corollary}
\newtheorem{rem}[teo]{Remark}
\newtheorem{ass}{Assumption}
\newtheorem*{defi}{Definition}

\DeclareFontFamily{OT1}{pzc}{}
\DeclareFontShape{OT1}{pzc}{m}{it}{<-> s * [1.10] pzcmi7t}{}
\DeclareMathAlphabet{\mathpzc}{OT1}{pzc}{m}{it}
\DeclareMathOperator{\Pic}{Pic}

\DeclareMathOperator{\id}{id}
\DeclareMathOperator{\discr}{discr}

\DeclareMathOperator{\Hom}{Hom}
\DeclareMathOperator{\divi}{div}

\DeclareMathOperator{\Mo}{Mo}

\newcommand{\pn}{\mathbb{P}}
\newcommand{\rk}{\mathrm{rk}\,}
\newcommand{\hsk}{K3^{\left[2\right]}-}
\newcommand{\IP}{ \mathbb{P}}
\newcommand{\IC }{ \mathbb{C}}
\newcommand{\IR}{ \mathbb{R}}

\newcommand{\IZ}{\mathbb{Z}}
\newcommand{\IQ}{\mathbb{Q}}
\newcommand{\IN}{\mathbb{N}}
\newcommand{\IH}{\mathbb{H}}

\newcommand{\prf}{{\it Proof.~}}
\newcommand{\ra}{\rightarrow}
\newcommand{\eprf}{\hfill $\square$ \smallskip\par}
\begin{document}
 \title{Lattice polarized irreducible holomorphic symplectic manifolds}
\author{Chiara Camere}
\address{Chiara Camere, Leibniz Universit\"at Hannover,
Institut f\"ur Algebraische Geometrie,
Welfengarten 1
30167 Hannover, Germany} 
\email{camere@math.uni-hannover.de}
\urladdr{http://www.iag.uni-hannover.de/~camere}
\date{\today}
\maketitle
\begin{abstract}
We generalize  Nikulin's and Dolgachev's lattice-theoretical mirror symmetry for K3 surfaces to lattice polarized higher dimensional irreducible holomorphic symplectic manifolds. In the case of fourfolds of $\hsk$type we then describe mirror families of polarized fourfolds and we give an example with mirror non-symplectic involutions.
\end{abstract}

\section{Introduction}

One striking prediction about geometrical objects coming from physics is the mirror conjecture.
Mirror symmetry for holomorphic symplectic manifolds has already been studied by Verbitsky in \cite{V}, where he shows that general non-projective holomorphic
symplectic manifolds are mirror self-dual; nothing is known  about projective holomorphic symplectic manifolds
apart for the two-dimensional case of K3 surfaces.

In \cite{DolgachevMirror} Dolgachev, based on former work by Nikulin \cite{Nikulin}, develops a mirror construction for lattice polarized projective K3 surfaces. First of all, he defines a moduli space $\mathfrak{M}_M$ parametrizing $M-$polarized K3 surfaces, i.e. those $S$ such that $M$ is
primitively embedded in $\mathrm{Pic}(S)$. Then he shows that, whenever there is a decomposition $M^{\perp}\cap H^2(S,\mathbb{Z})=U(m)\oplus \check{M}$, where $U$ is the standard hyperbolic lattice and $m$ is
an integer, then $\mathfrak{M}_{\check{M}}$ is a mirror moduli space: its dimension equals the Picard number of the generic member of $\mathfrak{M}_M$ and viceversa.
Moreover, the Griffiths-Yukawa coupling $\mathrm{Y}: S^2(H^1(S,T_S))\rightarrow H^{0,2}(S)^{\otimes 2}$ is a
symmetric pairing and for some open subset $\mathcal{U}$ of a compactification of $\mathfrak{M}_M$ near a boundary point, it can be identified with the quadratic form on
$\check{M}\otimes \mathbb{C}$. 
Finally, the period map of K3 surfaces induces a holomorphic multivalued map, the {\it mirror map}, from the open set $\mathcal{U}$ above to the tube domain
$\mathrm{Pic}(X')_{\mathbb{R}}+i\mathcal{K}_{X'}$, where $X'\in \mathfrak{M}_{\check{M}} $ and $\mathcal{K}_{X'}$ is its K\"ahler cone. 

Interesting examples of such a duality are given by Dolgachev, e.g. mirror partners of polarized K3s and Arnold's Strange Duality, and also by Borcea \cite{Borcea} and Voisin \cite{Voisin}, who introduced the notion of mirror non-symplectic involutions.

In this paper we generalize the definition of this lattice-theoretical mirror construction to higher dimensional irreducible holomorphic symplectic manifolds. After reviewing the basic notions of lattice theory and of the theory of hyperk\"ahler manifolds, in Section \ref{msslp} we construct moduli spaces of marked lattice polarized irreducible holomorphic symplectic manifolds and study their period domains. Then in Section \ref{mirror} we show how the theory in \cite{DolgachevMirror} carries through to higher dimensions: we define mirror moduli spaces so that they share the same properties mentioned above; roughly speaking, this duality exchanges the complex and the K\"ahler structure of the manifolds. In the case of fourfolds of $\hsk$type we then describe mirror families of polarized fourfolds and we generalize also the notion of mirror non-symplectic involutions.

{\it Acknowledgements.} The author wants to thank Klaus Hulek for suggesting this problem and for many enlightening discussions. She is also grateful to Samuel Boissi\`{e}re and Alessandra Sarti for their precious comments. This work was developed while the author was a member of the DFG Research Training Group ``Analysis, Geometry and String Theory'' at Leibniz University Hannover.

\section{Preliminary notions}
\subsection{Lattices}

A {\it lattice} $L$ is a free $\IZ$-module equipped with a nondegenerate symmetric bilinear form
$( \cdot, \cdot)$ with integer values.  Its {\it dual lattice} is 
$L^{*}:=\Hom_{\IZ}(L,\IZ)$ and can also be described as $
L^{*}\cong\{x\in L\otimes \IQ~|~( x,v)\in \IZ\quad \forall v\in L\}.
$
Since $L$ is a sublattice of $ L^{*}$ of the same rank, the quotient $A_L:=L^{*}/L$ is a finite abelian group, so-called {\it discriminant group}, of order $\discr(L)$, the {\it discriminant of $L$}. We denote by
$\ell(A_L)$ the minimal number of generators of $A_L$ (i.e. the \emph{length} of $A_L$).  In a basis $\{e_i\}_i$ of $L$, if $M:=(( e_i,e_j))_{i,j}$ is a Gram matrix, one has $\discr(L)=|\det(M)|$.

A lattice $L$ is called \emph{even} if $( x,x)\in 2\IZ$ for all $x\in L$.  In this case   
the bilinear form induces a finite quadratic form $q_L: A_L\longrightarrow \IQ/2\IZ$. If $(t_{(+)},t_{(-)})$ is the signature of
$L\otimes\IR$, the triple of invariants $(t_{(+)},t_{(-)},q_L)$ characterizes the \emph{genus} of the even lattice $L$ (see \cite[\S 7, Ch. 15]{conwaysloane}, 
\cite[Corollary 1.9.4]{Nikulinintegral}).

A lattice $L$ is called {\it unimodular} if $A_L=\{0\}$. An embedding of a sublattice $i:M\subset L$ is called \emph{primitive} if $L/i(M)$ is a free $\IZ$-module.
If $L$ is unimodular and $M\subset L$ is a primitive sublattice, then $M$ and its orthogonal $M^\perp$ in $L$ have
isomorphic discriminant groups and $q_M=-q_{M^\perp}$. When $L$ is no longer unimodular, the picture becomes more complicated and the following result helps with finding all non-isomorphic primitive embeddings of $M$.

\begin{teo}\cite[Proposition 1.15.1]{Nikulinintegral}\label{Nikiembeddings}. 
 The primitive embeddings of $M$ with invariants $(m_{(+)},m_{(-)},q_M)$ into an even lattice $L$ with invariants $(t_{(+)} ,t_{(-)}, q_L)$ are determined by the sets $(H_M, H_L, \gamma; K , \gamma_K)$ satisfying the following conditions:
\begin{itemize}
\item $H_M$ is a subgroup of $A_M$, $H_L$ is a subgroup of $A_L$ and $\gamma\colon H_M\to H_L$ is an isomorphism of groups such that for any $x\in H_M$, $q_L(\gamma(x))=q_M(x)$.
\item $K$ is a lattice of invariants $(t_{(+)}-m_{(+)},t_{(-)}-m_{(-)},q_K)$ with $q_K=\left.\left((-q_M)\oplus q_L\right)\right|_{\Gamma^\perp/\Gamma}$,
where $\Gamma$ is the graph of $\gamma$ in $A_M\oplus A_L$, $\Gamma^{\perp}$ is the orthogonal complement of $\Gamma$ in
$A_M\oplus A_L$ with respect to the bilinear form induced on $A_M\oplus A_L$ and with values in $\IQ/\IZ$; finally
$\gamma_K$ is an automorphism of $A_K$ that preserves $q_K$. Moreover $K$ is the orthogonal complement of $M$
in $L$. 
\end{itemize}
Two such sets, $(H_M, H_L, \gamma; K , \gamma_K)$ and $(H'_M, H'_L, \gamma'; K' , \gamma_{K'})$, determine isomorphic
primitive embeddings if and only if \begin{enumerate}
                                     \item $H_M=H'_M$;
                                     \item there exist $ \xi \in O(q_L)$ and $\psi:K\ra K'$ isomorphism for which $\gamma' = \xi \circ \gamma$ and $\xi\circ\gamma_K=\gamma_{K'}\circ\bar{\psi}$, where $\bar{\psi}$ is the isomorphism of the discriminant forms $q_K$ and $q_{K'}$ induced by $\psi$.
                                    \end{enumerate}

\end{teo}

In this paper $U$ will be the unique even unimodular hyperbolic 
lattice of rank two and $A_k, D_h, E_l$ will be the even, negative definite lattices associated to the Dynkin diagrams 
of the corresponding type ($k\geq 1$, $h\geq 4$, $l=6,7,8$). 
For $d\equiv -1\ (4)$ the following negative definite lattice will be used in the sequel
$$
K_d :=\left(
\begin{array}{cc}
-(d+1)/2&1\\
1&-2\\
\end{array}
\right) 
$$
Moreover, $L(t)$ denotes the lattice whose bilinear form is the one on $L$ multiplied by $t\in\IN^\ast$. 

We recall the following result by Nikulin on splitting of lattices.
\begin{teo}\cite[Theorem 1.13.5]{Nikulinintegral}\label{Niki}. Let $L$ be an even indefinite lattice 
of signature $(t_{(+)},t_{(-)})$ and assume that $t_{(+)}>0$ and $t_{(-)}>0$. Then:
\begin{enumerate}
\item If $t_{(+)}+t_{(-)}\geq 3+\ell(A_L)$ then $L\cong U\oplus W$ for a certain even lattice $W$.\\
\item If  $t_{(-)} \geq 8$ and $t_{(+)}+t_{(-)}\geq 9+\ell (A_L)$ then $L\cong E_8\oplus W'$ for a certain even lattice $W'$.
\end{enumerate}
\end{teo}

Finally, recall that the {\it divisor} $\divi f$ of a primitive element $f\in L$ is the generator of the ideal $(f,L)$ in $\IZ$.

\begin{teo}\cite[Prop. 3.7.3, Eichler's criterion]{Scattone}\label{eichler}
 If $L$ contains $U\oplus U$, given two primitive elements $f,f'\in L$ such that $f^2=(f')^2$ and $\divi f=\divi f'$ then there is an isometry $\sigma\in O(L)$ such that $\sigma(f)=f'$.
\end{teo}

\subsection{Irreducible holomorphic symplectic manifolds}

Irreducible holomorphic symplectic manifolds, also called hyperk\"ahler manifolds, have received a growing interest since it is known that if $X$ is a compact simply connected K\"ahler manifold with $c_1(X)_{\mathbb{R}}=0$ then there is a finite \'etale cover of $X$ that is a product of manifolds of three different types, namely complex tori, Calabi-Yau's and irreducible holomorphic symplectic ones (see \cite{B}). 

 A compact K\"{a}hler manifold $ X $ is {\it irreducible holomorphic symplectic} if it is simply connected and admits a symplectic two-form $ \omega_X \in
H^{2,0}(X) $ everywhere non degenerate and unique up to multiplication by a nonzero scalar.
The existence of such a symplectic form $\omega_X$ immediately implies that $ \dim X $ is an even integer. Moreover, $ K_X$ is trivial, in particular $c_1(X)=0$, and 
 $ T_X \cong \Omega _X^1 $. For a complete survey of this topic we refer the reader to the nice book \cite{GrossJoyceHuy} and references therein.
 
 The group $H^2(X,\mathbb{Z})$ carries a natural structure of lattice; the quadratic form on it is the so-called Beauville quadratic form $q$, which is even in all known examples.
We briefly recall here all known deformation types of irreducible holomorphic symplectic manifolds.
\begin{description}
\item[K3 surfaces]  These are compact complex connected surfaces $S$ with $b_1(S)=0$ and trivial canonical bundle. 
There is a lattice isomorphism between $H^2(S,\mathbb{Z})$ endowed with the cup-product and the lattice $3U\oplus 2E_8$.
\item[The Hilbert scheme of a K3 surface] Let $ S $ be a smooth K3 surface and let $ X=S^{\left[2\right]} $ be the Hilbert scheme of $ S $
of $
0- $schemes of length 2; $ X $ can be constructed also
as the blow-up along the image of the diagonal $ \Delta $ of the symmetric product $S^{(2)}$.
In particular, $ b_2(X)=23 $ and $h^{1,1}(X)=21$.  The Beauville lattice $(H^2(X,\mathbb{Z}),q)$ is $L=3U\oplus 2E_8\oplus \langle -2\rangle$.
Often, irreducible holomorphic symplectic manifolds that are deformation equivalent to $X$ are called \textit{of $K3^{\left[2\right]}-$ type}.

If $S$ is projective, then $X$ is also and $\Pic(X)\cong\Pic(S)\oplus\IZ e$ where $e$ is the class of the exceptional divisor with square $e^2=-2$.

The construction can be generalized in dimension $2n$, taking the Hilbert scheme of $S$ of $0-$schemes of length $n$ (see \cite{B}).

\item[Generalized Kummer varieties]  Given a complex torus $A$, we can construct the associated Hilbert scheme of $0-$schemes of length
$r+1$
and we get  $\xymatrix{A^{\left[r +1\right]}\ar@{->}[r]^{\ \ \ \pi} &A^{(r +1)} \ar@{->}[r]^{\ \ \ S}& A}$, where $S $ is the addition map. The
generalized Kummer variety is $ K^r(A):= (S\circ \pi)^{-1} (0)$ and it is an irreducible holomorphic symplectic manifold of dimension $2r$, as shown by Beauville in \cite{B}.
The Beauville lattice of a generalized Kummer variety is $3U\oplus\langle -2(r+1)\rangle$.

\item[O'Grady's examples] O'Grady in \cite{OG10} and \cite{OG6} discovered two deformation families of irreducible holomorphic symplectic manifolds of dimension
$6$ and $10$: these are deformations of desingularized moduli spaces of sheaves on abelian and K3 surfaces respectively.
\end{description} 

A \textit{marking} of an irreducible holomorphic symplectic manifold is an isometry $\phi:H^2(X,\IZ)\ra L$, where $L$ is a fixed even non-degenerate lattice of signature $(3,b_2(X)-3)$; a pair $(X,\phi)$ is then said to be {\it marked}. Similarly to what happens for K3 surfaces, there exists a moduli space of marked irreducible holomorphic symplectic manifolds $\mathcal{M}_L$ of fixed deformation type and one can define a \textit{period map} $\mathcal{P}_0:\mathcal{M}_L\ra D_L$ such that $\mathcal{P}_0(X,\phi)=\left[\phi(H^{2,0}(X))\right]$ in the \textit{period domain} $D_L:=\{\left[\omega \right]\in \pn(L_{\mathbb{C}})\mid (\omega,\omega)=0,(\omega,\bar{\omega})>0 \}$.

Already in \cite{B} it was shown that the period map is a local isomorphism; later Huybrechts in \cite{Huy1} show that $\mathcal{P}_0$ is surjective, even when restricted to a connected component. Finally, Verbitsky in \cite{VerTorelli} proved the global Torelli theorem, that here we recall as it has been formulated by Markman (see also Huybrechts's Bourbaki talk \cite{HuyBourbaki}).

\begin{teo}\cite[Theorem 1.3]{Mar}\label{HTT}. 
 Let $X$ and $Y$ be two irreducible holomorphic symplectic manifolds deformation equivalent one to each other. Then:
 \begin{enumerate}
  \item $X$ and $Y$ are bimeromorphic if and only if there exists a parallel transport operator $f:H^2(X,\IZ)\ra H^2(Y,\IZ)$ that is an isomorphism of integral Hodge structures;
  \item if this is the case, there exists an isomorphim $\tilde{f}:X\ra Y$ inducing $f$ if and only if $f$ preserves a K\"ahler class.
 \end{enumerate}
\end{teo} 

For the definition of parallel transport operator we refer to \cite[Definition 1.1]{Mar}, where also monodromy operators are defined.
\begin{defi}
 Given a marked pair $(X,\phi)$ of type $L$, we define the \textit{monodromy group} as $\Mo^2(L):=\phi\circ \Mo^2(X)\circ\phi^{-1}$, where $\Mo^2(X)\subset GL(H^2(X,\IZ))$ is the group of monodromy operators of $X$ restricted to the second cohomology group.
\end{defi}

It was proven by Verbitsky in \cite{VerTorelli} that $\Mo^2(L)$ is an arithmetic subgroup of $O(L)$; on the other hand we do not have an explicit description of this group in all known examples.

Moduli spaces of marked irreducible holomorphic symplectic manifolds are not Hausdorff, but we know exactly how to describe non-separated points.

\begin{teo}\cite[Theorem 2.2]{Mar}\label{GTT}.
 Let $\mathcal{M}^0_L$ be a fixed connected component of $\mathcal{M}_L$.
 \begin{enumerate}
  \item The period map $\mathcal{P}_0$ restricted to $\mathcal{M}^0_L$ is surjective.
  \item The fibre $\mathcal{P}_0^{-1}(p)$ consists of pairwise inseparable points for all $p\in D_L$.
   \item The marked pair $(X,\phi)$ is a Hausdorff point of $\mathcal{M}_L$ if and only if the positive cone and the K\"ahler cone coincide.  
 \end{enumerate}

\end{teo}

Things behave better when one restricts oneself to moduli of \textit{polarized} marked irreducible holomorphic symplectic manifolds, which have been studied in \cite{GHS} (see in particular Theorem 1.5). We denote with $\mathcal{A}_X$ the ample cone of $X$. Furthermore, given a primitive element $h\in L$, let $\Mo^2(h):=\left\lbrace g\in\Mo^2(L)\mid g(h)=h\right\rbrace$ be the subgroup of $h-$polarized monodromy operators and let $\Gamma_h$ be the image of $\Mo^2(h)$ via the restriction map $\alpha:O(L)\ra O(h^{\perp})$.

\begin{teo}\cite[Theorem 8.4]{Mar}\label{GHSPol}. Let $h\in L$ be a primitive element and let $D^+_h$ be one of the two connected components of $D_L\cap \IP(h^{\perp})$.
Let $\mathcal{M}_h^+$ be a connected component of the moduli space $\left\lbrace(X,\phi)\in\mathcal{P}_0^{-1}(D_h^+)\mid\phi^{-1}(h)\in\mathcal{A}_X\right\rbrace$ of polarized marked pairs. Then the period map restricts to an open embedding  with dense image$$\mathcal{P}_h:\mathcal{M}_h^+/\Mo^2(h)\ra D^+_h/\Gamma_h$$
\end{teo}

\section{Moduli spaces of lattice polarized irreducible holomorphic symplectic manifolds}\label{msslp}

This construction aims to generalize the one by Gritsenko, Hulek and Sankaran in \cite{GHS} for polarized irreducible holomorphic symplectic manifolds. Here we treat the subject in full generality and we will then specialize it to the case of fourfolds of $\hsk$type in Section \ref{hilb2}.

Let $X$ be an irreducible holomorphic symplectic manifold of type $L$ and let $j:M\subset L$ be a fixed primitive embedding of a sublattice $M$ of signature $(1,t)$; we will freely identify $M$ with $j(M)$ whenever confusion is not possible.

\begin{defi}
 An $M-$polarization of an irreducible holomorphic symplectic manifold $X$ is a lattice embedding $i:M\rightarrow \Pic(X)$. 
 
 A $j-$marking of an $M-$polarized manifold $X$ is a marking $\phi: H^2(X,\IZ)\rightarrow L$ such that $\phi\circ i=j$; a pair $(X,\phi)$ with $X$ an $M-$polarized irreducible holomorphic symplectic manifold of given deformation type and $\phi$ a $j-$marking is said to be $(M,j)-$polarized.

If $\Pic(X)=i(M)$ we say that $(X,\phi)$ is strictly $(M,j)-$polarized.
\end{defi}

Since $i(M_{\mathbb{C}})\subset H^{1,1}(X)$, we have that $\mathcal{P}_0(X,\phi)\in\pn (M_{\mathbb{C}}^{\perp})$; hence, we can consider a restricted
period domain $$D_{M}=\{\left[\omega \right]\in \pn(N_{\mathbb{C}})\mid
(\omega,\omega)=0,(\omega,\bar{\omega})>0 \}$$ where $N=j(M)^{\perp}$. This has two connected components and each one is a symmetric homogenous domain of type IV (see \cite{GHSHandbook}).
 Since $N_{\IC}$ depends only on the signature of $N$, the period domain $D_M$ depends only on $M$ and not on $j$.

We need a notion of ample polarization and to introduce it we need to make the following assumption:
\begin{ass}\label{ass1}
There exists a set $\Delta(L)\subset L$ such that the K\"ahler cone $\mathcal{K}_X$ of a marked $(X,\phi)$ can be described as $$\mathcal{K}_X=\left\lbrace h\in H^{1,1}(X,\mathbb{R})\mid (h,h)>0, (h,\delta)>0\  \forall\delta\in\Delta(X)^+\right\rbrace$$
where $\Delta(X)^+:=\left\lbrace\delta\in \phi^{-1}(\Delta(L))\cap\Pic(X)\mid (\delta,\kappa)>0\right\rbrace$ for $\kappa\in \mathcal{K}_X$ a fixed K\"ahler class.
\end{ass}
As we will see more in detail in Section \ref{hilb2}, Assumption \ref{ass1} is satisfied in the case of fourfolds of $K3^{\left[2\right]}-$type.
Given such an embedding $j:M\subset L$, define the positive cone $C(M)=\left\lbrace x\in M_{\mathbb{R}}\mid (x,x)>0\right\rbrace$ and pick one of the two connected components $C^+(M)$. Given $\Delta(M):=\Delta(L)\cap M$ and $H_{\delta}=\left\lbrace x\in M_{\mathbb{R}}\mid (x,\delta)=0 \right\rbrace$, we fix a connected component of $C^+(M)\setminus(\cup_{\delta\in\Delta(M)}H_{\delta})$ and call it $K(M)$. This choice induces the choice of a set $\Delta(M)^+=\left\lbrace\delta\in\Delta(M)\mid (x,\delta)>0 \ \forall x\in K(M)\right\rbrace$ such that $\Delta(M)=\Delta(M)^+\amalg(-\Delta(M)^+)$.

\begin{defi}
 We say that  $(X,\phi)$ as above is ample (strictly) $(M,j)-$polarized if $i(K(M))$ contains a K\"ahler class.
\end{defi}

\begin{lem}\label{ampgenMpol}
If $(X,\phi)$ is ample strictly $(M,j)-$polarized then:
\begin{enumerate}
 \item $i(\Delta(M)^+)=\Delta(X)^+$;
 \item $i(K(M))=\mathcal{K}_X$.
\end{enumerate}
\end{lem}
\prf (1) Take $\kappa\in i(K(M))\cap\mathcal{K}_X$; given $\delta\in i(\Delta(M)^+) $, we have $(\delta,\kappa)=$ $(\tilde{\delta},\tilde{k})>0$ for $\tilde{\delta}=\phi(\delta)\in \Delta(M)^+ $ and $\tilde{k}=\phi(\kappa)\in K(M)$. On the other hand, suppose that there exists $\delta\in\Delta(X)^+\setminus i(\Delta(M)^+)$; we have that $\delta\in i(\Delta(M)^-)$ and hence $(\delta,\kappa)=(\tilde{\delta},\tilde{k})<0$, in contradiction with our assumption.

(2)This follows immediately from (1) and the definitions.\eprf

Given a smooth family $f:\mathcal{X}\rightarrow \mathcal{U}$ of irreducible holomorphic symplectic manifolds of given deformation type, an $M-$polarization of $f$ is an injection $i_\mathcal{U}:M_{\mathcal{U}}\rightarrow \mathcal{P}\mathpzc{ic}_{\mathcal{X}/S}\subset R^2f_*\IZ$, from the constant sheaf $M_{\mathcal{U}}$ to the relative Picard sheaf $\mathcal{P}\mathpzc{ic}_{\mathcal{X}/S}$, such that for every $t\in \mathcal{U}$ the map $i_t$ defines an $M-$polarization of $\mathcal{X}_t$. A $j-$marking of the family is then defined (see \cite{DolgachevMirror}) as an isomorphism of local systems $\phi_{\mathcal{U}}:R^2f_*\IZ\rightarrow L_{\mathcal{U}}$ such that $\phi_t\circ i_t=j$ for all $t\in {\mathcal{U}}$. Such a marking allows to define the {\it period map} of $f$ as $\mathcal{P}_f:t\in {\mathcal{U}}\mapsto \left[\phi_t(\omega_{\mathcal{X}_t})\right]\in D_M$, that is holomorphic by the local Torelli Theorem \cite[Th\'eor\`eme 5]{B}. Local moduli spaces and period maps are then glued together and give a coarse moduli space $\mathcal{M}_{M,j} $ of $(M,j)-$polarized irreducible holomorphic symplectic manifolds of fixed type and a holomorphic map $\mathcal{P}_{M,j}:\mathcal{M}_{M,j}\rightarrow D_M$ that is the restriction of the period map $\mathcal{P}_0$ of marked irreducible holomorphic symplectic manifolds of fixed type.

The group $O(L,M)=\{g\in O(L)\mid g(m)=m\ \forall m\in M\}$ acts properly and discontinuously on $D_M$; choose a connected component $D^+_{M}$ of $D_M$ and take  
$\mathcal{M}_{M,j}^+$ a connected component of $\mathcal{P}_{M,j}^{-1}(D^+_{M})$; it is a connected component of $\mathcal{M}_{M,j}$ and the period map restricts to a surjective holomorphic map $\mathcal{P}_{M,j}^+:\mathcal{M}_{M,j}^+\rightarrow D_{M}^+$ which is a local isomorphism.

As defined in Markman \cite{Mar}, consider $\Mo^2(L):=\phi\circ \Mo^2(X)\circ \phi^{-1}$ where $(X,\phi)\in \mathcal{M}_L^+$. We define {\it $(M,j)-$polarized monodromy operators} $$\Mo^2(M,j):=\{g\in \Mo^2(L)\mid g(m)=m\ \forall m\in M\}=\Mo^2(L)\cap O(L,M) $$ 
In other words, an element $g\in\Mo^2(M,j)$ satisfies $g\circ j=j$. This group acts on $\mathcal{M}_{M,j}^+$ via $(X,\phi)\mapsto (X,g\circ\phi)$ for $g\in\Mo^2(M,j)$. 

The restriction map induces an injective map $\alpha:g\in O(L,M)\mapsto g_{|N}\in O(N)$; we define the subgroup $\Gamma_{M,j}:=\alpha(\Mo^2(M,j))$ in $O(N)$.

\begin{pro} \label{Mpol}
The set $\mathcal{M}_{M,j}^+$ is invariant under the action of $\Mo^2(M,j)$ and the restriction of the period map is $\Mo^2(M,j)-$equivariant, so that we get a surjective map
\[
 \xymatrix{
\mathcal{M}_{M,j}^+/\Mo^2(M,j)\ar[r]^{\mathcal{P}_{M,j}^+} & D^+ _{M}/\Gamma_{M,j}
}
\]
\end{pro}

\prf Given $(X,\phi)\in \mathcal{M}_{M,j}^+$ and $g\in\Mo^2(M,j)$, there is an embedding $i:M\subset\Pic(X)$ such that $\phi\circ i=j$; then $g\circ\phi$ is again a $j-$marking since $g\circ\phi\circ i=g\circ j=j$.

The equivariance of the restricted period map is trivial.\eprf

To obtain a quasi-projective variety we need to show that $\Gamma_{M,j}$ is of finite index inside $O(N)$. 
By a result of Markman combined with work of Kneser (see also \cite{GHSHandbook}), it follows that if $X$ is of $K3^{\left[2\right]}-$type then $\Mo^2(L)$ is related to the so-called {\it stable orthogonal group}, $$\tilde{O}^+(L)=\left\lbrace g\in O(L)\mid g_{|A_L}=\id, \mathrm{sn}_{\IR}^L(g)=1\right\rbrace$$
where the {\it real spinor norm} $\mathrm{sn}_{\IR}^L:O( L_{\IR})\rightarrow \IR^*/(\IR^*)^2\cong\lbrace \pm 1\rbrace$ is defined as $$\mathrm{sn}_{\IR}^L(g)=\left(-\frac{v_1^2}{2}\right)\cdots\left(-\frac{v_m^2}{2}\right)(\IR^*)^2$$ for $g\in O(L_{\IR})$ factored as a product of reflections $g=\rho_{v_1}\circ\dots\circ\rho_{v_m}$ with $v_i\in L_{\IR}$.

\begin{pro}\label{fi}
If $\Mo^2(L)\supset \tilde{O}^+(L)$ the group $\Gamma_{M,j}$ is an arithmetic subgroup of $O(N)$.
\end{pro}
\prf We will show that $H=\tilde{O}^+(N)$ is a subgroup of $\Gamma_{M,j}$ and since it is of finite index in $O(N)$ (this follows immediately from the fact that $\mathrm{Aut}(A_N)$ is a finite group), $\Gamma_{M,j}$ is also of finite index in $O(N)$ because $O(N)/H\supseteq O(N)/\Gamma_{M,j}$.

Given $g\in \tilde{O}^+(N)$ we want to prove that there exists $f\in \Mo^2(M)$ such that $\alpha(f)=g$. Take $f\in O(L)$ to be the map induced on $L$ by $\id_M\oplus g$; then by definition $f\in O(L,M)$. Moreover, $f_{|(A_M\oplus A_N)}=\id_{A_M\oplus A_N}$, since $g\in \tilde{O}(N)$, and $A_L\subset A_M\oplus A_N$ (from $M\oplus N\subset L\subset L^*\subset M^*\oplus N^*$), hence $f_{|A_L}=\id_{A_L}$ and $f\in\tilde{O}(L)$.

Next, consider the extension of $g$ by linearity to $N_{\IR}$; we know that there are $v_1,\dots,v_m\in N_{\IR}$ such that $g=\rho_{v_1}\circ\dots\circ\rho_{v_m}$ in $O(N_{\IR})$ and $\mathrm{sn}_{\IR}^N(g)=1$. We will still denote by $\rho_{v_i}$ the reflection of $L_{\IR}$ with respect to $v_i\in N_{\IR}\subset L_{\IR}$; for all $v_i$ with $i=1,\dots,m$ we have that $(\rho_{v_i})_{|M_{\IR}}=\id_{M_{\IR}}$ since $(v_i,m)=0$ for all $m\in M_{\IR}$, hence also $f=\rho_{v_1}\circ\dots\circ\rho_{v_m}$ in $O(L_{\IR})$ and $\mathrm{sn}_{\IR}^L(f)=\mathrm{sn}_{\IR}^N(g)=1$, i.e. $f\in O^+(L)$.

So indeed, $f\in  \Mo^2(M,j)$ and by construction $\alpha(f)=g$.\eprf

\begin{cor}
If $\Mo^2(L)\supset \tilde{O}^+(L)$ the quotient $D_{M,j}^+/\Gamma_{M,j}$ is a quasi-projective variety of dimension $\rk L-2-\rk M$. 
\end{cor}

\prf This follows from Proposition \ref{fi} and Baily-Borel's theorem \cite{BailyBorel}.\eprf

Next we restrict to a connected component $\mathcal{M}_{M,j}^a\subset \mathcal{M}_{M,j}^+$ of the moduli space of ample $(M,j)-$polarized irreducible holomorphic symplectic manifolds and we denote by $\mathcal{M}_{M,j}^{sa}$ the subset of ample strictly $(M,j)-$polarized ones.

\begin{lem}
The set $ \mathcal{M}_{M,j}^a$ is open in the analytic topology.
\end{lem}
\prf We have that $\mathcal{M}_{M,j}^a=\cup (\mathcal{M}_{M,j}^a\cap \mathcal{M}_{h})$ and once a polarization $h$ is fixed, ampleness is an open condition due to to the stability of K\"ahler manifolds (see \cite[\S9.3.3]{VoisinBook}).\eprf
\begin{teo}\label{ampgenMpolmod}
 The subset $\mathcal{M}_{M,j}^{sa}$ is Hausdorff in $\mathcal{M}_{M,j}^+$ and it is invariant under the action of $\Mo^2(M,j)$. Moreover, its image via the period map is connected and dense.
\end{teo}
\prf First we show that $$\mathcal{P}_{M,j}^+(\mathcal{M}_{M,j}^{sa})= D^{\circ}_{M}:= D_{M}^+\setminus \left(\bigcup_{\nu\in N\setminus\left\lbrace0 \right\rbrace} H_{\nu}\right),$$
where $H_{\nu}=\left\lbrace \lambda\in N_{\IC}\mid(\lambda,\nu)=0\right\rbrace$.

Indeed, given $(X,\phi)\in \mathcal{M}_{M,j}^{sa}$ and $\pi=\mathcal{P}_{M,j}^+(X,\phi)\in D_{M}$, we see that $\pi\notin H_{\nu}$ for any $\nu\in N\setminus\left\lbrace0 \right\rbrace$: otherwise, $\nu\in\pi^{\perp}$ and $\phi^{-1}(\nu)\in \Pic(X)\setminus i(M)$, contradicting our assumption.

Given $\pi\in  D^{\circ}_{M}$ and $(X,\phi)\in\mathcal{P}_{M,j}^{-1}(\pi)$, the pair $(X,\phi)$ is strictly $(M,j)-$polarized; moreover, there is a bijection $\rho:\mathcal{P}_{M,j}^{-1}(\pi)\rightarrow\mathcal{KT}(X)$ via $\rho(Y,\eta)=\eta^{-1}(\phi(\mathcal{K}_X))$  by \cite[Proposition 5.14]{Mar}, where $\mathcal{KT}(X)$ is the set of K\"ahler-type chambers of $X$. If $i(K(M))\cap\mathcal{K}_X\neq\emptyset$ there is nothing to prove; otherwise, $i(K(M))$ meets a different K\"ahler-type chamber since $\Delta(X)=i(\Delta(M))$. Hence, there exists $(Y,\eta)\in \mathcal{M}_{M,j}^{sa}\cap \mathcal{P}_{M,j}^{-1}(\pi)$; in fact, it follows easily from Theorem \ref{GTT} and Lemma \ref{ampgenMpol} that there exists a unique such $(Y,\eta)$, so that $\mathcal{M}_{M,j}^{sa}$ is Hausdorff.

Finally remark that $D^{\circ}_{M}$ is connected and dense by Baire's category theorem.\eprf

\begin{cor}
 The period map induces a bijection 
 \begin{displaymath}
 \mathcal{P}_{M,j}^{sa}: \mathcal{M}_{M,j}^{sa}/\Mo^2(M,j)\longrightarrow D^{\circ}_{M}/\Gamma_{M,j}                                      
 \end{displaymath}

\end{cor}
\prf It follows from Proposition \ref{Mpol} and Theorem  \ref{ampgenMpolmod} that the period map restricts to a bijection $\mathcal{P}_{M,j}^{sa}: \mathcal{M}_{M,j}^{sa}\rightarrow D^{\circ}_{M} $ and that the restriction is equivariant with respect to the action of $\Mo^2(M,j)$. \eprf

\begin{rem}
 If the primitive embedding $j:M\subset L$ is unique up to isometry of $L$, then $\mathcal{M}_M$ can be seen as the moduli space of $M-$polarized
irreducible holomorphic symplectic manifolds, getting rid of markings as done by Dolgachev in \cite{DolgachevMirror}. On the other hand, since $L$ is no longer unimodular, this is a stronger
condition to require with respect to the case of K3 surfaces as it is not always satisfied even in the case of polarizations (see \cite{GHS}). Proposition \ref{uniqembed} describes some cases in which this happens.
\end{rem}

\section{Mirror symmetry}\label{mirror}
\subsection{Griffiths-Yukawa coupling}

In this section we limit ourselves to recalling some notations and facts from \S4 in \cite{DolgachevMirror} and we focus our attention on the few modifications needed in higher dimensions.

From now on suppose that $\rk M\leq 20$, so that its orthogonal $N$ in $L$ (which is unique up to isometry once we fix the embedding $j:M\hookrightarrow L$ by Theorem \ref{Nikiembeddings}) is indefinite. Fix a primitive isotropic vector $f\in N_{\IR}$, so that $(f,f)=0$, and set 
$$
N_f=\left\lbrace x\in N_{\IR}\mid (x,f)=1\right\rbrace,\ V_f=\left\lbrace x\in N_{\IR}\mid (x,f)=0\right\rbrace/\IR f;                                                                                                                                                                                                                                                                                                                                                                                                                                                                                                          $$
let $C^+_f$ be a connected component of the cone
$$
C_f=\left\lbrace x\in V_f\mid (x,x)>0\right\rbrace
$$
The corresponding \textit{tube domain}, which is the complexification of $C^+_f$, is 
$$\mathcal{H}_f=N_f+iC^+_f $$

\begin{pro}\cite[Corollary 4.3]{DolgachevMirror}
The choice of an  isotropic $f\in N_{\IR}$ determines an isomorphism $D^+_{M}\cong \mathcal{H}_f$. 
\end{pro}

Now we want to use such a tube domain realization to relate the Griffiths-Yukawa coupling with the quadratic form on $N$. We are dealing with manifolds of dimension $2n$, but the period domain $D^+_{M}$ parametrizes weight-two Hodge structures and as such we are still interested in looking at the Griffiths-Yukawa quadratic form 
\begin{displaymath}
 Y:S^{2}H^1(X,T_X)_{\phi}\longrightarrow H^{0,2}(X)^{\otimes 2}\end{displaymath}
 $$
 (\theta_1,\theta_2)\mapsto \varphi^{1,1}(\theta_1)\circ\varphi^{2,0}(\theta_2)
$$
where $H^1(X,T_X)_{\phi}$ is the tangent space of $\mathcal{M}_{M,j}^+$ at the point $(X,\phi)$, defined as the orthogonal in $H^1(X,T_X)$ of $i(M)$ with respect to the pairing
$$
 H^1(X,T_X)\otimes H^{1,1}(X)\longrightarrow H^{0,2}(X),
$$
and $\varphi^{i,j}:H^1(X,T_X)\rightarrow \Hom (H^{i,j}(X),H^{i-1,j+1}(X))$, for $1\leq i\leq 2n$ and $0\leq j\leq 2n-1$, is given by the interior product with a tangent vector (see \cite{Griffiths}).

\begin{pro}(\cite[Corollary 4.4]{DolgachevMirror})
 For any $\mu\in D^+_{M}$ the choice of a representative $l\in L$ of $\mu$ such that $(l,f)=1$ defines a canonical isomorphism $$\alpha_{\mu}:T_{\mu}D_{M}^+\rightarrow (V_f)_{ \IC}$$
 
 Moreover, given $(X,\phi)\in \mathcal{M}_{M,j}^+$, the quadratic form on $(V_f)_{\IC}$ coincides with the Griffiths-Yukawa pairing with respect to the normalization $H^{0,2}(X)\cong\IC$ defined by $\phi^{-1}(l)\in H^{2,0}(X)$.
\end{pro}

\subsection{The mirror map}

Again the theoretical construction contained in \S5 and \S6 of \cite{DolgachevMirror} carries over to higher dimensions with very little modification. First of all, let us recall some definitions.
\begin{defi}
 Let $S$ be an even indefinite lattice and $m$ a positive integer; a primitive isotropic vector $f\in S$ is $m-$admissible if $\divi f=m$ and there exists another isotropic vector $g\in S$ such that $(f,g)=m$, $\divi g=m$.
\end{defi}
Due to \cite[Lemma 5.4]{DolgachevMirror} this is equivalent to the existence of a primitive embedding $U(m)\rightarrow S$ such that $f\in U(m)$.

Suppose that there exists an $m-$admissible $f\in N$; then $N=U(m)\oplus \check{M}$, with $\check{M}$ a primitive sublattice of signature $(1, h^{1,1}-\rk M-1)$, and we have $\check{M}\cong (\IZ f)^{\perp}_{N}/\IZ f$. The definition of $\check{M}$ depends only on the choice of $f$ if $U(m)$ admits a unique primitive embedding in $N$: this happens in particular if $(m,\det N)=1$, $\check{M}$ is unique in its genus and there is a surjection $O(\check{M})\rightarrow O(q_{\check{M}})$. Once this holds, we can define $\check{\jmath}$ to be the composition of the embedding $\check{M}\subset N$ and of the embedding $j^\perp:N\subset L$ and this depends only on $f$ and $j$.

\begin{defi}
The moduli space $\mathcal{M}_{\check{M},\check{\jmath}}^+$ is the mirror moduli space of $\mathcal{M}_{M,j}^+$. 
\end{defi}

\begin{pro}
We have that $\dim \mathcal{M}_{\check{M},\check{\jmath}}^+=\rk M$ and $\dim \mathcal{M}_{\check{M},\check{\jmath}}^++\dim \mathcal{M}_{M,j}^+=h^{1,1}$.
\end{pro}

We introduce now the Baily-Borel compactification of the period domain, defined as its closure in the Harish-Chandra embedding (see \cite{GHSHandbook} for a nice survey of the topic), contained in the obvious compactification given by the quadric $$D^*_{M}=\{\left[\omega \right]\in \pn(N_{\mathbb{C}})\mid
(\omega,\omega)=0\}.$$ A {\it boundary component} is a subset of the form $\IP(I_{\IC})\cap D_{M}^*$ for some isotropic subspace $I\subset N_{\IR}$ of dimension $1$ or $2$; such a component is called \textit{rational} if the corresponding $I$ can be defined over $\IQ$. In particular, $0-$dimensional rational boundary components of $D^+_{M}$ are in bijection with primitive isotropic elements of $N$.

When $\Gamma_{M,j}$ is an arithmetic subgroup of $O(N_{\IQ})$, it acts on the set of rational boundary components of $D^+_{M}$, which we denote $\mathcal{RB}$, and for each such $F$ its stabilizer $N(F)=\{g\in \Gamma_{M,j}\mid g(F)=F\}$ acts discretely on $F$. Then the Baily-Borel compactification is
\begin{displaymath}
 \overline{D^+_{M}/\Gamma_{M,j}}=D^+_{M}/\Gamma_{M,j}\coprod \left(\coprod_{F\in\mathcal{RB}/\Gamma_{M,j}}F/N(F)\right)
\end{displaymath}
that is a normal projective algebraic variety.

Now we choose an $m-$admissible primitive isotropic $f\in N$ and consequently we fix a splitting $N=U(m)\oplus \check{M}$ and an isotropic $g\in U(m)$ such that $(f,g)=m$. We define $Z_{M,j}(f)=\{h\in N(F)\mid h(f)=f\}$ and $Z_{M,j}(f)^+$ as the subgroup of elements preserving $K(\check{M})$; thus we have an action of $Z_{M,j}(f)$ on $\mathcal{H}_f=V_f+iC_f^+$ and we can identify $Z_{M,j}(f)^+$ with the subgroup preserving $\mathcal{H}^+_f=V_f+iK(\check{M})$.

Let $F$ be the $0-$dimensional rational boundary component corresponding to $f$. The theory in \cite{DolgachevMirror} holds also in our situation, hence there exist neighbourhoods $\tilde{{\mathcal{U}}}^*$ and ${\mathcal{U}}$ respectively of $F$ in $D^*_{M}$ and of $F/N(F)$ in $\overline{D^+_{M}/\Gamma_{M,j}}$ such that
$\alpha: \tilde{{\mathcal{U}}}/Z_{M,j}(f)^+\rightarrow {\mathcal{U}}_F\subset {\mathcal{U}}$ is an isomorphism for $\tilde{{\mathcal{U}}}=\tilde{{\mathcal{U}}}^*\cap\mathcal{H}_f^+$.

\begin{teo}
 The period map induces a multi-valued map, called \textit{mirror map}, $\alpha^{-1}: {\mathcal{U}}_F\rightarrow \tilde{{\mathcal{U}}}_F\subset \mathcal{H}_f^+$ with monodromy group equal to $Z_{M,j}(f)^+$ that sends a neighbourhood of $F$ to the tube domain $\mathcal{H}_f^+\cong \Pic(X')+i\mathcal{K}_{X'}$ for $(X',\phi')\in\mathcal{M}^{sa}_{\check{M},\check{\jmath}}$.
\end{teo}
\section{The $K3^{\left[2\right]}-$type case}\label{hilb2}

From now on let $X$ be a fourfold of $K3^{\left[2\right]}-$type, so that $b_2(X)=23$ and $L=3U\oplus 2E_8\oplus \left\langle-2\right\rangle$. In this case, by works of Hassett and Tschinkel \cite{HassTschink} and of Mongardi \cite{MongardiMori} (see also \cite{BayerHassettTschinkel}), Assumption \ref{ass1} is satisfied with $$\Delta(L)=\left\lbrace\delta\in L\mid \delta^2=-2\right\rbrace\cup\left\lbrace\delta\in L\mid \delta^2=-10, \divi(\delta)=2\right\rbrace.$$

By a result of Markman combined with work of Kneser (see also \cite{GHSHandbook}), it follows that  $$\Mo^2(L)=\mathrm{Ref}(L)=\tilde{O}^+(L)=\left\lbrace g\in O(L)\mid g_{|A_L}=\id, \mathrm{sn}_{\IR}^L(g)=1\right\rbrace$$
Hence the hypothesis of Proposition \ref{fi} is satisfied and $\Gamma_{M,j}$ is an arithmetic subgroup of $O(N)$. 

\begin{teo}\label{fik32}
 For $L=3U\oplus 2E_8\oplus \left\langle-2\right\rangle$ we have that $\Gamma_{M,j}=\tilde{O}^+(N)$.
\end{teo}
\prf Proposition \ref{fi} tells us that $\Gamma_{M,j}\supset\tilde{O}^+(N)$. Viceversa, in this case $\Mo^2(M,j)=\tilde{O}^+(L)\cap O(L,M)$. Given $f\in\Mo^2(M,j)$, then $f_{|A_M\oplus A_L}=id$ and hence $f_{|A_N}=id$, since by Nikulin's Theorem \ref{Nikiembeddings} $\Gamma, \Gamma^{\perp}\subset A_M\oplus A_L$. On the other hand, $\Gamma_{M,j}\subset O^+(N)$ because it preserves the connected component $D^+_{M}$, hence we have equality.\eprf

In this case we can find some criteria for the unicity of the embedding $j$.
\begin{pro}\label{uniqembed}
If there is no subgroup $H\subset A_M$ such that $H\cong\IZ/2\IZ$, the orthogonal $N$ is unique in its genus $(2,20-t,q_N)$ and the projection $O(N)\rightarrow O(q_N)$ is surjective, then $M$ admits a unique primitive embedding $j:M\hookrightarrow L$. 
\end{pro}

In particular this is true if $M$ is unimodular of rank $\rk M\leq 20$ or if $A_M=\bigoplus_{p_i>2\mathrm{prime}} (\IZ/p_i\IZ)^{\oplus a_i}$ and $\rk M\leq 21-\max(a_i)$. It is important to stress though that the orthogonal $N$ will not in general have a unique embedding, so that $\check{\jmath}$ will not be the only possible embedding of $\check{M}$.

\begin{rem}
 Consider now a primitive embedding $j_{K3}:M\subset L_{K3}$ and take $j=j_{K3}\oplus\id_{\left\langle -2\right\rangle}$ to be the induced primitive embedding $M\subset L$. 
We can then find a mirror lattice either in $L_{K3}$ or in $L$, obtaining respectively $\check{M}_{K3}$ and $\check{M}=\check{M}_{K3}\oplus \left\langle -2\right\rangle$. Take an $M-$polarized K3 surface $S$ and $S'$ an $\check{M}_{K3}-$polarized K3 surface in the mirror family; then $S^{\left[2\right]}$ and $(S')^{\left[2\right]}$ will be respectively $M-$polarized and $\check{M}-$polarized mirror partners. On the other hand, since for any $M-$polarized K3 surface $S$, $S^{\left[2\right]}$ is also $(M\oplus \left\langle -2\right\rangle)-$polarized, the family of Hilbert schemes of $M-$polarized K3 surfaces has codimension $1$ inside $\mathcal{M}^+_{M,j}$, whereas the mirror moduli space has the same dimension as in the K3 case.
\end{rem}

\subsection{The polarized case}

Let $M\subset L$ be the rank $1$ sublattice $\left\langle2d\right\rangle$ for $d$ a positive integer. In the $K3^{\left[2\right]}-$type case, the following result is known:

\begin{teo}\cite[Prop. 3.6 and 3.12]{GHS}
 The sublattice $M=\left\langle2d\right\rangle$ admits up to two non-isometric primitive embeddings in $L$. Let $h$ be a generator of $M$; then the following hold:
 \begin{enumerate}
  \item there is always the {\it split} embedding $j_s$, corresponding to $\divi h=1$, such that $N_s=2U\oplus2E_8\oplus \left\langle-2\right\rangle\oplus \left\langle-2d\right\rangle$, $\det N_s=4d$ and $A_{N_s}=\IZ/2\IZ\oplus \IZ/2d\IZ$;
  \item if $d\equiv 3$ modulo $4$ then $M$ admits a second embedding  $j_{ns}$, called {\it non-split}, corresponding to $\divi h=2$, such that $N_{ns}=2U\oplus2E_8\oplus K_d$, $\det N_{ns}=d$ and $A_{N_{ns}}=\IZ/d\IZ$.
 \end{enumerate}
In both cases, $\Gamma_{M,j}\cong \tilde{O}(N)^+$.
\end{teo}

Fix $j:M\subset L$ and let $\mathcal{M}_{2d}^+$ be the subset of $\mathcal{M}^+_{M,j}$ where $h$ corresponds to an ample class; as already shown in \cite{GHS}, we get an isomorphism with a Zariski open set.

\begin{pro}
 The period map $\mathcal{P}_{M,j}^+$ restricts to an isomorphim $$\mathcal{M}_{2d}^+/\Mo^2(M,j)\rightarrow \left(D^+_{M}\setminus\coprod_{\delta\in\Delta(N)}(H_{\delta}\cap D^+_{M})\right)/\tilde{O}(N)^+ $$
\end{pro}

Now we want to compute $\tilde{O}(N)-$orbits of $m-$admissible isotropic vectors $f$ in $N$ for an integer $m|\det N$. In both cases we can apply Eichler's criterion \ref{eichler}: orbits are classified by $\divi f|\det N$. From Scattone's work\cite[Proposition 4.1.3]{Scattone} there is a bijection between $\tilde{O}(N)$-orbits of isotropic vectors in $N$ and the set of isotropic elements in $A_N$ modulo multiplication by $\pm 1$, induced by the map $f\in N\mapsto f/\divi f+N\in A_N$.

{\bf The split case}. Let $e$ and $t$ denote respectively the generators of $\left\langle-2\right\rangle$ and $\left\langle-2d\right\rangle$ in $N_{s}$. Then the discriminant group $A_{N_{s}}$ is generated by $e/2$ and $t/2d$ and the discriminant quadratic form $q_{s}$ is given by 
\begin{displaymath}
q_{s}(\alpha \frac{t}{2d}+\beta \frac{e}{2})=-\frac{\alpha^2+\beta^2d}{2d}\in\IQ/2\IZ                                                                                                                                                                                                                                                                                                                                                                                                                                                                                 \end{displaymath}
for $\alpha=0,\dots,2d-1$ and $\beta=0,1$. 
Let $u$ and $v$ denote a standard basis of one of the two orthogonal summands $U$ inside $N_s$ and write $d=d'k^2$ with $d'$ square-free.
\begin{lem}\label{misot}
The $m-$admissible isotropic vectors $f$ in $N_s$, up to the action of $\tilde{O}(N_s)$, are of the form
\begin{equation}\label{isot}
 f=\left\lbrace
 \begin{array}{cc}
t+m(u+\frac{d}{m^2}v)&\mathrm{if}\ m|k\\
&\\
t+m^2e+m(u+\frac{d+m^2}{m^2}v)&\mathrm{if}\ m\nmid k,\ \frac{m}{2}|k\ \mathrm{and}\ d'\equiv 3 (4)
\end{array}\right.
\end{equation}
\end{lem}
\prf Let $I(q_s)$ be the set of isotropic elements in $A_{N_s}$; it is clear that it is the union of isotropic elements of $\IZ/2d\IZ$ with respect to the restricted form $q_s(\alpha \frac{t}{2d})=-\frac{\alpha^2}{2d}$ and of isotropic elements $y=\alpha \frac{t}{2d}+ \frac{e}{2}$, since $e$ is not isotropic. Moreover, it is easy to remark that  $\mathrm{ord}(f/\divi f)=\divi f$ and, since we are interested in classification up to the action of $\tilde{O}(N_s)$, we need to find only one isotropic element for each possible order $m$.

The computation of isotropic elements of $\IZ/2d\IZ$ has been done by Scattone \cite[Theorem 4.0.1]{Scattone} in the K3 case; they are in bijection with elements in the cyclic subgroup of order $k$. Hence $m$ has to divide $k$ and an isotropic element $x\in\IZ/2d\IZ$ of order $m $ is $x_m=\frac{t}{m}$.

Now take $y=\alpha \frac{t}{2d}+ \frac{e}{2}$; then $q_s(y)=-\frac{\alpha^2+d}{2d}\in 2\IZ$ if and only if $\alpha^2+d=4dl$ for an integer $l$. By reducing this equality modulo $4$ we see that the only possible case is $d'\equiv 3\ (4)$ and $\alpha=d'kh$ with $h=1,\dots, k-1$ odd.
In this case we observe that $m=\mathrm{ord}(y)|2k$ and we are interested in finding elements of order not dividing $k$; in particular, new possible orders $m$ are those such that $\frac{m}{2}|k$ and $m\nmid k$. For example, an isotropic element $x\in A_{N_s}$ of order $m $ is $x_m=\frac{t}{m}+\frac{e}{2}$.

Up to isometry, the corresponding isotropic vectors in $N_s$ are precisely of the forms given in (\ref{isot}). Indeed, given such an $f\in N_s$, we have $f^2=0$, $(f,v)=m$,
\begin{displaymath}
 (f,u)=\left\lbrace
 \begin{array}{cc}
\frac{d}{m}\in m\IZ&\mathrm{if}\ m|k\\
&\\
\frac{d+m^2}{m}\in m\IZ&\mathrm{if}\ m\nmid k,\ \frac{m}{2}|k\ \mathrm{and}\ d'\equiv 3 (4)
\end{array}\right.,
\end{displaymath}
\begin{displaymath}
  (f,e)=\left\lbrace
 \begin{array}{cc}
0&\mathrm{if}\ m|k\\
&\\
-2m^2\in m\IZ&\mathrm{if}\ m\nmid k,\ \frac{m}{2}|k\ \mathrm{and}\ d'\equiv 3 (4)
\end{array}\right.
\end{displaymath}
Hence $\divi f=m$.

\eprf

Now we restrict to the case $m=1$ and consider the isotropic primitive vector $f=t+u+dv$. Since it is unimodular, $U$ admits a unique primitive embedding in $N$ up to isometry; hence its orthogonal is $\check{M}=U\oplus 2E_8\oplus \left\langle-2\right\rangle\oplus \left\langle-2d\right\rangle$ and we can assume that $U=\IZ f+\IZ v$. The sublattice $\check{M}$ admits two non-isometric primitive embeddings in $L$; it follows from the definition that $\check{\jmath}$ satisfies $\check{\jmath}(\check{M})^{\perp}= U\oplus \left\langle2d\right\rangle$. Hence the period domain 
$D^+_{\check{M}}$ is exactly the one described by Dolgachev in \cite{DolgachevMirror}, with tube domain realization isomorphic to the upper half-plane $\IH$.

From \cite[Theorem 7.1]{DolgachevMirror} and Theorem \ref{fik32}, the global monodromy group $\Gamma_{\check{M},\check{\jmath}}$ is conjugate in $PSL(2,\IR)$ to $ \Gamma_0(d)^+$ generated by the group  $$\Gamma_0(d)=\left\lbrace (a_{ij})\in SL(2,\IZ)\mid d|a_{21}\right\rbrace$$ and by the Fricke involution 
\begin{displaymath}                                                                                                                                                                                                                                                                                      F=\left(\begin{array}{cc}                                                                                                                                                                                                                                                                                        0&-\frac{1}{\sqrt{d}}\\                                                                                                                                                                                                                                                                                        \sqrt{d}&0                                                                                                                                                                                                                                                                                                                                                                                                                                                                                                                                                                                      \end{array}
\right)\in PSL(2,\IR)                                                                                                                                                                                                                                                                                        \end{displaymath}

The main difference with what happens in the case of polarized K3 surfaces is that here we only get a local isomorphism from our moduli space to the modular curve $$ \xymatrix{\mathcal{M}_{\check{M},\check{\jmath}}^+/\Mo^2(\check{M},\check{\jmath})\ar[r]^{\mathcal{P}} & \IH/\Gamma_0(d)^+}$$

{\bf The non-split case.} In this case $d\equiv 3\ (4)$. Let $e$ and $w_1$, $w_2$ denote respectively the generators of $\left\langle-2\right\rangle$ and of a copy of $U$ in $L$ so that $M=\left\langle h \right\rangle\subset U\oplus\left\langle -2\right\rangle$ via $h=2w_1+\frac{d+1}{2}w_2+e$ with orthogonal $K_d$ generated by $b_1=w_1-\frac{d+1}{4}w_2$ and $b_2=w_2+e$. The discriminant group $A_{N_{ns}}$ is generated by $t=\frac{1}{d}h-w_2$ and the discriminant quadratic form $q_{ns}$ is given by 
\begin{displaymath}
q_{ns}(\alpha t)=-\frac{2\alpha^2}{d}\in\IQ/2\IZ                                                                                                                                                                                                                                                                                                                                                                                                                                                                                 \end{displaymath}
for $\alpha=0,\dots,d-1$. 
Let $u$ and $v$ denote a standard basis of one of the two orthogonal summands $U$ inside $N_{ns}$ and write $d=d'k^2$ with $d'\equiv 3\ (4)$ square-free.
\begin{lem}
The $m-$admissible isotropic vectors $f$ in $N_{ns}$, up to the action of $\tilde{O}(N_{ns})$, are of the form
\begin{equation}\label{nsisot}
 f=2b_1+b_2+m(u+\frac{d}{m^2}v)\ \mathrm{if}\ m|k
\end{equation}
\end{lem}
\prf Computations similar to the ones in the proof of Lemma \ref{misot} show that the order $m$ of an isotropic element $\alpha t\in \IZ/d\IZ$ has to divide $k$. Given $f$ as in (\ref{nsisot}), we have that $f^2=0$, $(f,b_1)=-d\in m\IZ$, $(f,b_2)=0$, $(f,u)=\frac{d}{m}\in m\IZ$ and $(f,v)=m$. Hence $\divi f=m$.\eprf

Now we restrict to the case $m=1$ and consider the hyperbolic lattice $U=\IZ f\oplus\IZ v$. Since it is unimodular, $U$ admits a unique primitive embedding in $N$ up to isometry and its orthogonal is $\check{M}=U\oplus 2E_8\oplus K_d$. The sublattice $\check{M}$ admits a unique primitive embedding into $L$ and $\check{\jmath}(\check{M})^{\perp}= U\oplus \left\langle2d\right\rangle$; the period domain $D^+_{\check{M}}$ is exactly as in the split case and everything remarked above holds again.

\subsection{Non-symplectic involutions} 

In the forthcoming paper \cite{BCS} the authors classify primitive embeddings of invariant sublattices $T$ of non-symplectic involutions $i$ of fourfolds $X$ of $K3^{\left[2\right]}-$type, i.e. involutions such that $i^*\omega_X=-\omega_X$. The invariant sublattice $T$ is known to be hyperbolic and two-elementary with two-elementary orthogonal $S$. By work of Nikulin \cite{Nikulinfactor} a two-elementary hyperbolic lattice $T$ is completely determined by the triple $(r,a_T,\delta_T)$, where $r$ is its rank, $a_T=l(A_T)$ is the length of its discriminant group and $\delta_T$ is the {\it parity} of the discriminant quadratic form $q_T$: $\delta_T=0$ if $q_T(x)\in\IZ/2\IZ$ for all $x\in A_T$, $1$ otherwise.

By \cite[Proposition 6.1]{BCS} primitive embeddings $j$ into $L$ of a two-elementary hyperbolic sublattice $T$ with invariants $(r,a_T,\delta_T)$ are in bijection with couples $(a\pm 1,\delta_S)$ where $S$ is the orthogonal complement of $j(T)$ in $L$, two-elementary with $a_S=l(A_S)=a\pm 1$ and parity $\delta_S$. 

Consider now $\mathcal{M}^+_{T,j}$ and look for the mirror family corresponding to the choice of a $1-$admissible isotropic vector $f\in S$.

\begin{lem}
There is a $1-$admissible isotropic $f\in S$ and $S\supset U$ if and only if $\rk T\leq 21-l(A_T)$ or $\rk T=22-l(A_T)$ and $l(A_S)=l(A_T)+1$ except for $(15,7,1,6,0)$.
\end{lem}
\prf If $\rk T\leq 20-l(A_S)$ then this follows by Theorem \ref{Niki}. Otherwise, one of the following hold:
\begin{enumerate}
 \item $\rk T=24-l(A_T)$ and $l(A_S)=l(A_T)-1=\rk S$: we have that $S=\left\langle 2\right\rangle^{\oplus 2}\oplus \left\langle -2\right\rangle^{\oplus l(A_T)-3}$ if $\delta_S=1$ and if $\delta_S=0$ then $S$ is $2U(2)$. None contains a copy of $U$.
 \item  $\rk T=22-l(A_T)$ and $l(A_S)=l(A_T)-1=\rk S-2$: we have that $S=U\oplus \left\langle 2\right\rangle\oplus \left\langle -2\right\rangle^{\oplus l(A_T)-2}$ if $\delta_S=1$ and if $\delta_S=0$ then $S$ is either $U\oplus U(2)$,  $2U(2)\oplus D_4$ or $U\oplus U(2)\oplus E_8(2)$. All contain a copy of $U$ except $2U(2)\oplus D_4$.
 \item  $\rk T=22-l(A_T)$ and $l(A_S)=l(A_T)+1=\rk S$: we have that $S=\left\langle 2\right\rangle^{\oplus 2}\oplus \left\langle -2\right\rangle^{\oplus l(A_T)-1}$ with $\delta_S=1$. None contains a copy of $U$.
 \item  $\rk T=20-l(A_T)$ and $l(A_S)=l(A_T)+1=\rk S-2$: we have that $S=U\oplus\left\langle 2\right\rangle\oplus \left\langle -2\right\rangle^{\oplus l(A_T)}$ with $\delta_S=1$ and this contains a copy of $U$.
\end{enumerate}\eprf

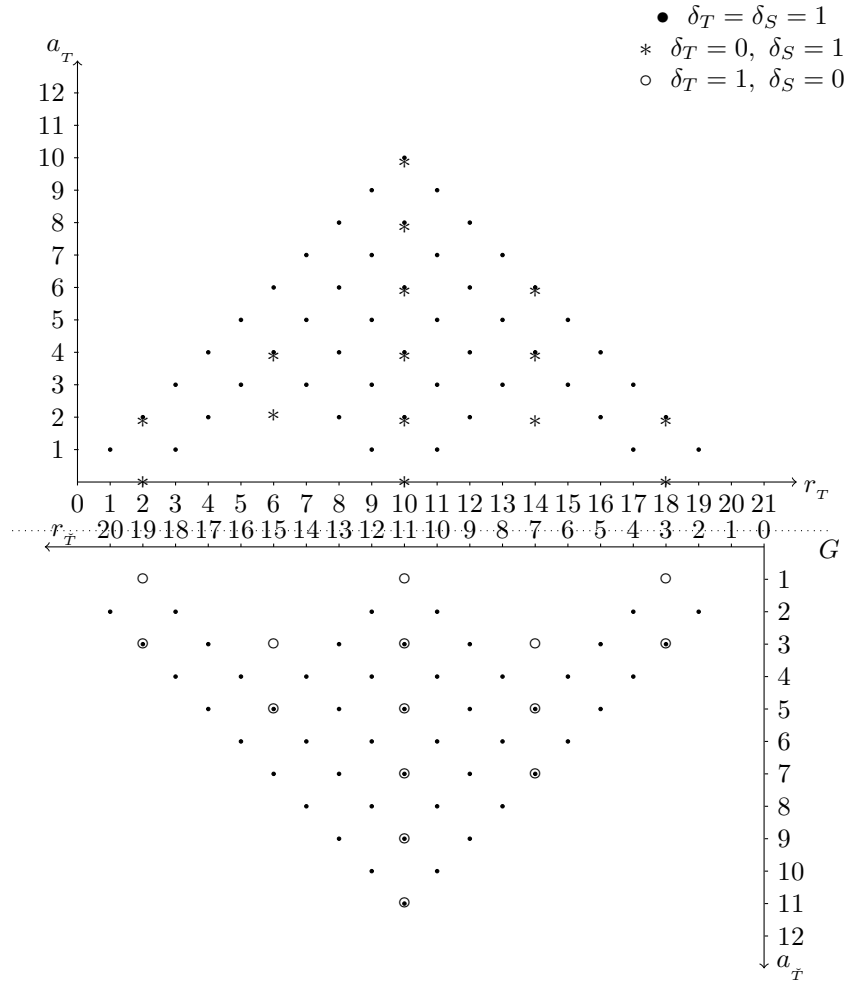
\begin{figure}[h]
$$\begin{array}{cccccccccccccccccccccccr}
 &\ &\ &&&& &&\ &\ &&&&&&&&&&&&& \bullet\ \ \delta_T=\delta_S=1\\
&\ &\ &&&& &&\ &\ &&&&&&&&&&&&&\ast\ \ \delta_T=0,\ \delta_S=1\\
&\ &\ &&&& &&\ &\ &&&&&&&&&&&&&\circ\ \ \delta_T=1,\ \delta_S=0
\end{array}$$
\vspace{-1.1cm}

\begin{tikzpicture}[scale=.43]
\filldraw [black] 
(1,1) circle (1.5pt)  
(2,0) node[below=-0.20cm]{*} 
(2,2) circle (1.5pt) node[below=-0.15cm]{*}
 (3,1) circle (1.5pt)
 (3,3) circle (1.5pt)
 (4,2) circle (1.5pt)
(4,4) circle (1.5pt)
(5,3) circle (1.5pt)
(5,5) circle (1.5pt)
(6,4) circle (1.5pt)node[below=-0.15cm]{*}
(6,2)    node[below=-0.23cm]{*}
(6,6) circle (1.5pt)
(7,3) circle (1.5pt)
(7,5) circle (1.5pt)
(7,7) circle (1.5pt)
(8,2) circle (1.5pt)
(8,4) circle (1.5pt)
(8,6) circle (1.5pt)
(8,8) circle (1.5pt)
(9,1) circle (1.5pt)
(9,3) circle (1.5pt)
(9,5) circle (1.5pt)
(9,7) circle (1.5pt)
(9,9) circle (1.5pt)
(10,0)  node[below=-0.20cm]{*}
(10,2) circle (1.5pt)node[below=-0.15cm]{*}
(10,4) circle (1.5pt)node[below=-0.15cm]{*}
(10,6) circle (1.5pt)node[below=-0.15cm]{*}
(10,8) circle (1.5pt)node[below=-0.15cm]{*}
(10,10) circle (1.5pt)node[below=-0.15cm]{*}
(11,1) circle (1.5pt)
(11,3) circle (1.5pt)
(11,5) circle (1.5pt)
(11,7) circle (1.5pt)
(11,9) circle (1.5pt)
(11,11) 
(12,2) circle (1.5pt)
(12,4) circle (1.5pt)
(12,6) circle (1.5pt)
(12,8) circle (1.5pt)
(12,10)
(13,3) circle (1.5pt)
(13,5) circle (1.5pt)
(13,7) circle (1.5pt)
(13,9) 
(14,2)  node[below=-0.15cm]{*}
(14,4) circle (1.5pt)node[below=-0.15cm]{*}
(14,6) circle (1.5pt)node[below=-0.15cm]{*}
(14,8)
(15,3) circle (1.5pt)
(15,5) circle (1.5pt)
(15,7) 
(16,2) circle (1.5pt)
(16,4) circle (1.5pt)
(16,6) 
(17,1) circle (1.5pt)
(17,3) circle (1.5pt)
(17,5) 
(18,0)  node[below=-0.2cm]{*}
(18,2) circle (1.5pt)node[below=-0.15cm]{*}
(18,4) 
(19,1) circle (1.5pt)
(19,3) 
(20,2) 
 (19,-4) circle (1.5pt)
 (18,-3)node{$\circ$} 
 (18,-5) circle (1.5pt)node{$\circ$} 
 (17,-4) circle (1.5pt)
(17,-6) circle (1.5pt)
(16,-5) circle (1.5pt)
(16,-7) circle (1.5pt)
(15,-6) circle (1.5pt)
(15,-8) circle (1.5pt)
(14,-5) node{$\circ$} 
(14,-7) circle (1.5pt)node{$\circ$} 
(14,-9) circle (1.5pt)node{$\circ$} 
(13,-6) circle (1.5pt)
(13,-8) circle (1.5pt)
(13,-10) circle (1.5pt)
(12,-5) circle (1.5pt)
(12,-7) circle (1.5pt)
(12,-9) circle (1.5pt)
(12,-11) circle (1.5pt)
(11,-4) circle (1.5pt)
(11,-6) circle (1.5pt)
(11,-8) circle (1.5pt)
(11,-10) circle (1.5pt)
(11,-12) circle (1.5pt)
(10,-3) node{$\circ$} 
(10,-5) circle (1.5pt)node{$\circ$} 
(10,-7) circle (1.5pt)node{$\circ$} 
(10,-9) circle (1.5pt)node{$\circ$} 
(10,-11) circle (1.5pt)node{$\circ$} 
(10,-13) circle (1.5pt)node{$\circ$}  
(9,-4) circle (1.5pt)
(9,-6) circle (1.5pt)
(9,-8) circle (1.5pt)
(9,-10) circle (1.5pt)
(9,-12) circle (1.5pt)
(8,-5) circle (1.5pt)
(8,-7) circle (1.5pt)
(8,-9) circle (1.5pt)
(8,-11) circle (1.5pt) 
(7,-6) circle (1.5pt)
(7,-8) circle (1.5pt)
(7,-10) circle (1.5pt) 
(6,-5) node{$\circ$} 
(6,-7) circle (1.5pt)node{$\circ$} 
(6,-9) circle (1.5pt)
(5,-6) circle (1.5pt)
(5,-8) circle (1.5pt)
(4,-5) circle (1.5pt)
(4,-7) circle (1.5pt)
(3,-4) circle (1.5pt)
(3,-6) circle (1.5pt)
(2,-3) node{$\circ$} 
(2,-5) circle (1.5pt)node{$\circ$}  
(1,-4) circle (1.5pt)
 ; 
\draw[->] (0,0) -- coordinate (x axis mid) (22,0);
    \draw[->] (0,0) -- coordinate (y axis mid)(0,13);
    \foreach \x in {0,1,2,3,4,5,6,7,8,9,10,11,12,13,14,15,16,17,18,19,20,21}
        \draw [xshift=0cm](\x cm,0pt) -- (\x cm,-3pt)
         node[anchor=north] {$\x$};
          \foreach \y in {1,2,3,4,5,6,7,8,9,10,11,12}
        \draw (1pt,\y cm) -- (-3pt,\y cm) node[anchor=east] {$\y$};
    \node[below=0.1cm, right=4.7cm] at (x axis mid) {$r_{_T}$};
   \node[left=0.2cm, below=-3.2cm,] at (y axis mid) {$a_{_{_T}}$};
   \draw[->] (21,-2) -- coordinate (z axis mid) (-1,-2);
    \draw[->] (21,-2) -- coordinate (t axis mid)(21,-15);
    \foreach \x in {0,1,2,3,4,5,6,7,8,9,10,11,12,13,14,15,16,17,18,19,20}
        \draw [xshift=1cm](20 cm-\x cm,-54pt) -- (20cm-\x cm,-57pt) node[anchor=south] {$\x$};
          \foreach \y in {1,2,3,4,5,6,7,8,9,10,11,12}
        \draw (598pt,-\y cm-2cm) -- (601pt,-\y cm-2cm) node[anchor=west] {$\y$};
    \node[below=0.7cm, right=-5.2cm] at (x axis mid) {$r_{_{\check{T}}}$};
   \node[right=9.4cm, below=8.95cm,] at (y axis mid) {$a_{_{_{\check{T}}}}$};
  \draw[dotted,thin](-2,-1.5)--(23,-1.5)node[below, black]{$G$};
  \end{tikzpicture} 
\caption{Mirror pairs of non-symplectic involutions}\label{mir}
\end{figure}

Once fixed such an $f\in S$ and a splitting $S=U\oplus\check{T}$, we get that $\check{T}$ is hyperbolic, two-elementary with invariants $(21-r, a_S,\delta_S)$ and the embedding $\check{\jmath}$ is the one corresponding to $(a_T,\delta_T)$. Moreover, by cancelling the points corresponding to non-admissible values of $(r,a_T,\delta_T,a_S,\delta_S)$, Figure 1 and Figure 2 in \cite{BCS} can be combined in Figure \ref{mir}, where every point denoted with $\bullet$ is mirror dual with the symmetric $\bullet$ with respect to the line $G$ and symmetric $\ast$ and $\circ$ are mirrors.

\begin{rem}
Since the generic member of each family carries a non-symplectic involution with prescribed invariant lattice, we thus get a notion of {\it mirror involution}. On the other hand, this does not agree with the notion of mirror involutions defined by the analogous construction on K3 surfaces, as described in \cite{Voisin}, in the sense that pairs of natural involutions induced by mirror involutions on a K3 surface $S$ are not mirror pairs on $S^{\left[2\right]}$.
\end{rem}


\bibliographystyle{amsplain}
\bibliography{BiblioMSLP}

\end{document}